\documentclass[12pt]{article}
 \pdfoutput=1
  \usepackage{amssymb,amsmath,latexsym,theorem,eepic,enumerate}
  \usepackage[pagebackref=true,
              pdfpagelabels]{hyperref}
  \usepackage{times}
 \usepackage[all]{xy}  \usepackage{geometry,color}
{\theorembodyfont{\rmfamily}\newtheorem{example}{Example}[section]}
{\theorembodyfont{\rmfamily}\newtheorem{Def}[example]{Definition}}
{\theorembodyfont{\rmfamily}}
{\theorembodyfont{\rmfamily}}
{\theorembodyfont{\rmfamily}}

\newtheorem{theorem}[example]{Theorem}

\newcommand{\threeaxes}[3]{\def\objectstyle{\scriptstyle}  \objectmargin={0pt}
\xy
(0,0)*+{}="a",(0,-6)*+{\rule{0em}{1.5ex}#2}="b",(7,0)*+{\;#1}="c",
(14,-3)*+{\;#3}="d" \ar@{->} "a";"b" \ar @{->}"a";"c"  \ar
@{->}"a";"d"\endxy }

\newcommand{\directs}[2]{\def\objectstyle{\scriptstyle}  \objectmargin={0pt}
\xy
(0,4)*+{}="a",(0,-2)*+{\rule{0em}{1.5ex}#2}="b",(7,4)*+{\;#1}="c"
\ar@{->} "a";"b" \ar @{->}"a";"c" \endxy }

\newcommand{\xdirects}[2]{\def\objectstyle{\scriptstyle}  \objectmargin={0pt}
\xy
(0,0)*+{}="a",(0,-6)*+{\rule{0em}{1.5ex}#2}="b",(7,0)*+{\;#1}="c"
\ar@{->} "a";"b" \ar @{->}"a";"c" \endxy }
\newcommand{\sdirects}[2]{\def\objectstyle{\scriptstyle}  \objectmargin={0pt}
\xy
(0,2.2)*+{}="a",(0,-2.5)*+{\rule{0em}{1.5ex}#2}="b",(7,2.2)*+{\;#1}="c"
\ar@{->} "a";"b" \ar @{->}"a";"c" \endxy }

\newcommand{\bl}{\mbox{\rule{0.08em}{1.7ex}\hspace{-0.00em}\rule{0.7em}{0.2ex}}}

\newcommand{\br}{\mbox{\rule{0.7em}{0.2ex}\hspace{-0.04em}\rule{0.08em}{1.7ex}}}

\newcommand{\tr}{\mbox{\rule[1.5ex]{0.7em}{0.2ex}\hspace{-0.03em}\rule{0.08em}{1.7ex}}}

\newcommand{\tl}{\mbox{\rule{0.08em}{1.7ex}\rule[1.54ex]{0.7em}{0.2ex}}}

\newcommand{\hh}{\mbox{\rule{0.7em}{0.2ex}\hspace{-0.7em}\rule[1.5ex]{0.70em}{0.2ex}}}

\newcommand{\vv}{\mbox{\rule{0.08em}{1.7ex}\hspace{0.6em}\rule{0.08em}{1.7ex}}}

\newcommand{\sq}{\mbox{\rule{0.08em}{1.7ex}\hspace{-0.00em}\rule{0.7em}{0.2ex}\hspace{-0.7em}\rule[1.54ex]{0.7em}{0.2ex}\hspace{-0.03em}\rule{0.08em}{1.7ex}}}

\newcommand{\tsq}{\mbox{\rule{0.04em}{1.55ex}\hspace{-0.00em}\rule{0.7em}{0.1ex}\hspace{-0.7em}\rule[1.5ex]{0.7em}{0.1ex}\hspace{-0.03em}\rule{0.04em}{1.55ex}}}

\def\rho{\varrho}

\def\epsilon{\varepsilon}

\def\cal{\mathcal}

\def\w{\mathbf{w}}

\def\epsilon{\varepsilon}

\def\cal{\mathcal}

\def\bu{\bullet}
\def\cal{\mathcal}
\def\epsilon{\varepsilon}

\def\bu{\centerdot}

\def\xybiglabels{\def\labelstyle{\textstyle}}
\begin{document}

Mathematical Intelligencer

(old Submission id:  TMIN-D-21-00076  : revised) 

Subject: Manuscipt TMIN-D-20-0001OR1

entitled: ``Not just an idle game"\\(the story  of higher dimensional versions of the Poincar{\' e} fundamental group)

author: Ronald Brown

acknowledged by TMIN: March 6, 2021

\newpage
\setcounter{page}{1}
\title{``Not just an idle game"\\(the story  of higher dimensional versions of the Poincar{\' e}  fundamental group)}
\author{Ronald Brown}
\maketitle

\section*{Introduction}
Part of the  title of this article is taken from writings of
Einstein\footnote{More information on most of  the people mentioned in this article may be found in the web site \url{https://mathshistory.st-andrews.ac.uk/Biographies/}.} (1879-1955) in the correspondence published in \cite{Ein90}:

\begin{quote}
$\ldots$  the following questions must burningly interest me as a disciple of science: What goal will be reached by the science to which I am dedicating myself? What is essential and what is based only on the accidents of development? $\ldots$ Concepts which have proved useful for ordering things easily assume so great an authority over us, that we forget their terrestrial origin and accept them as unalterable facts.  $\ldots$  It is therefore not just an idle game to exercise our ability to analyse familiar concepts, and to demonstrate the conditions on which their justification and usefulness depend, and the way in which these developed, little by little $\ldots$
\end{quote}

This quotation is about science, rather than mathematics, and it is well known for example in physics that there are still fundamental questions, such as the nature of dark matter, to answer. There should be
awareness  in mathematics that there are still some basic questions which have failed to be pursued for decades; thus we need to think also of educational methods of encouraging their pursuit.

In particular,  encouragements or discouragements are personal, and may turn out not to be correct. That is a standard hazard of research. But it was the interest  of making reality of the intuitions of ``algebraic inverse to subdivision'' (see diagram \eqref{eq:aitos}) and of ``commutative cube'' (see Figure  \eqref{fig:commcub}) that kept the project  alive.

\section{Homotopy  groups at the ICM Z{\" u}rich, 1932}
Why am I considering this ancient meeting?  Surely we have advanced since then?   And the basic ideas have surely long been totally sorted?

Many mathematicians, especially Alexander\footnote{I make the convention that the use of a full first name indicates I had the good fortune to have had a significant   direct contact  with  the person  named. } Grothendieck (1928-2014),   have  shown us that basic ideas  can be looked at again and in some cases renewed.

The main  theme with which I am concerned in this paper  is little discussed today, but is stated in \cite[p.98]{jam}: it involves the introduction  by  the  well respected topologist E. {\v C}ech (1893-1960) of   {\it  homotopy groups}  $\pi_n(X,x)$ of a pointed space $(X,x)$, and which were proved by him  to be abelian  for $n > 1$.

But it was argued that these groups  were inappropriate for what was a key theme at the time,  the development of higher dimensional versions of the fundamental group $\pi_1(X,x)$ of a pointed space as defined by H. Poincar{\' e} (1854-1912), \cite{Gray}.  In many of the known applications of the fundamental group in complex analysis and differential equations, the largely nonabelian nature of the fundamental group was a key factor.

Because of this  abelian property  of higher homotopy groups, {\v C}ech    was persuaded by H. Hopf  (1894-1971)  to withdraw his paper, so that only a small paragraph appeared in the Proceeding \cite{Ce}.   However, the abelian homology groups $H_n(X)$ were known at the time to be well  defined for any space $X$,  and that if $X$ was path connected,  then $H_1(X)$ was isomorphic to  the fundamental group $\pi_1(X,x)$ made abelian.

Indeed P.S.  Aleksandrov (1896-1998)  was reported to have exclaimed: "But my dear \v Cech, how can they be anything but the homology groups?" \footnote{I heard of this comment  in Tbilisi in 1987 from  George  Chogoshvili  (1914-1998) whose doctoral supervisor was Aleksandrov. Compare also \cite{Alex}. An irony is that the 1931  paper \cite{Hopf} already gave a counterexample to this statement, describing   what is now known as the ``Hopf map''  $S^3 \to S^2$; this map is non trivial homotopically, but it is trivial homologically; it   also has very interesting relations to other aspects of algebra and geometry, cf. \cite[Ch. 20]{jam}.   Aleksandrov and Hopf were two of the most respected topologists: their standing is shown by their  mentions in \cite{jam} and the invitation by  S. Lefschetz (1884-1972) for them to spend the academic year 1926 in Princeton. }

{\bf Remark} It should be  useful to give here for $n=2$  an  intuitive argument which {\v C}ech might  have used for the abelian nature of the homotopy groups $\pi_n(X,x), n>1$.

{ \bf Proof} It is possible to represent any $g \in G=\pi_2(X,x)$ by a map $a: I^2 \to X$ which is constant with value $x$ outside a ``small'' square, say $J$, contained in $I^2$.

Another  such class $h \in G$ may be similarly  represented by a map $b$ constant  outside a small square $K$ in $I^2$, and such that $J$ does not meet $K$.  Now we can see a clear difference between the cases $ n=2, n=1$. In the former case, we can choose $J,K$ small enough and separated  so that the maps $a,b$ may be deformed in their classes so that  $J,K$ are interchanged.

 This method is not  possible if $n=1$, since in that case $I$ is an interval, and one of $J,K$  is to the left or right of the other. Thus it can seem that any expectations of higher dimensional versions of the fundamental group were unrealisable. Nowadays, this argument would be put in the form of what is called the Eckmann-Hilton Interchange result, that in a set with two monoid structures, each of which is a morphism for the other, the two structures coincide and are abelian.

In 1968,  Eldon Dyer  (1934-1997), a topologist at CUNY,   told me that Hopf told him in 1966 that the history of homotopy theory showed the danger of people being regarded as ``the kings" of a subject and so key in deciding directions. There is a lot in this point, cf \cite{Alex}.

Also at this ICM was W. Hurewicz (1906-1956); his publication  of  two notes \cite{Hur}  shed light on the relation of the homotopy groups to homology groups, and the interest in these homotopy groups started. With the growing study of the complications of the homotopy groups of spheres\footnote{This can be seen by a web search on this topic.}, which became seen as a major problem in algebraic topology, cf \cite{WG},  the idea of generalisations of the nonabelian fundamental group became disregarded, and it became easier to think of ``space" and the ``space with base point" necessary to define the homotopy groups, as in substance synonymous - that was my experience up to 1965.

 However it can be argued that Aleksandrov and Hopf {\it were  correct } in suggesting that  the abelian homotopy groups are not what one would really like for a higher dimensional generalisation of the fundamental group! That does not mean that such  ``higher homotopy groups''  would be  without interest; nor does it mean that the search for nonabelian higher dimensional generalisations of the fundamental groups should be completely abandoned.
\section{Determining fundamental groups}
One reason for this interest in fundamental groups was their known use in  important questions  relating  complex analysis, covering spaces, integration and group theory, \cite{Gray}. H.Seifert (1907-1996) proved useful relations between simplicial complexes and fundamental groups,  \cite{Sei}, and a paper by E.R.  Van Kampen (1908-1942), \cite{VK},  stated  a general result which could be applied to the complement in a 3-manifold of an algebraic curve.  A  modern proof for the case of path connected intersections was given  by Richard H. Crowell (1928-2006) in \cite{Cro} following lectures of R.H. Fox (1913-1973). That  result is often called the Van Kampen Theorem (VKT) and there are many excellent examples of applications of it in expositions  of algebraic topology.

The usual statement of the VKT for a fundamental group is as follows:
\begin{theorem}\label{thm:vktgp}
  Let the space $X$ be the union of open sets $U,V$ with intersection $W$ and assume $U,V,W$ are path connected. Let $x \in W$. Then the following diagram of fundamental groups and morphisms induced by inclusions:
  \begin{equation}\label{eq:vktgp}
 \vcenter{\xybiglabels \xymatrix{\pi_1(W,s) \ar_i [d] \ar^j [r] & \pi_1(V,x) \ar^k [d] \\
  \pi_1(U,x) \ar_h [r] & \pi_1(X,x), }}
  \end{equation}
   is a pushout diagram of groups.
\end{theorem}
Note that a  ``pushout of groups'' is defined entirely in terms of the notion of morphisms of groups: using  the  diagram \eqref{eq:vktgp} the definition  says that if $G$ is any group and  $f:\pi_1(U,x) \to G, g: \pi_1(V,x) \to G$ are morphisms of groups such that $fi=gj$, then there is a unique morphism of groups $\phi: \pi_1(X,x) \to G$ such that $\phi h= f, \phi g=k$.  This property is called the ``universal property'' of a pushout, and proving it is called ``verifying the universal property''.  It is often convenient that  such a verification need not involve a particular construction of the pushout, nor a proof that all pushouts of morphisms of groups exist. See also \cite{BJP}.

The limitation to path connected spaces and intersections in Theorem \ref{thm:vktgp} is also very restrictive.\footnote{Comments  by Grothendieck on this restriction, using the word ``obstinate'',  are quoted extensively in \cite{survey}, see also \cite[Section 2]{GrEsq}.}. Because of the connectivity condition on $W$, this  standard version of the Van Kampen Theorem for the fundamental group of a pointed space does not compute the fundamental group of the circle\footnote{The argument for this is that if $S^1=U \cup V$ is the union of two open, path connected sets, then $U \cap V$ has at least two path components.}, which is after all {\bf the}  basic example in topology; the standard treatments instead make a detour into a small part  of covering space theory by introducing the ``winding number'' of the map $p: \mathbb R \to S^1, t \mapsto e^{2\pi i t}$ from the reals to the circle, which goes back to Poincar\'e in the 1890s.

\section{From groups to groupoids}\label{sec:gpds}
This is a theme with which I became involved in the years since 1965.

 A {\it groupoid} is defined in modern terms as a small category in which every morphism is an isomorphism. It can be considered    as a ``group     with many identities'', or more formally as an algebraic structure with {\it partial}  algebraic operations, as considered  by Philip J. Higgins
 (1924-2015)  in \cite{Higgins1}.  I like to define ``higher dimensional algebra" as the study of partial algebraic structures where the  domains of the algebraic operations are defined by geometric conditions.

 Groupoids  had been defined by Brandt (1886-1954) in 1926, \cite{Brandt},  for  extending to the quaternary case work of Gauss (1777-1815) on compositions of binary quadratic forms;  the use of groupoids in topology had been initiated by K. Reidemeister (1893-1971)  in his 1932 book, \cite{Reid}.

 The simplest non trivial example of a groupoid is the groupoid say $\mathcal I$ which has two objects $0,1$ and only one nontrivial arrow $\iota:0 \to 1$, and hence also $\iota^{-1}: 1 \to 0$.  This groupoid looks ``trivial'', but it is in fact the basic ``transition operator''.  (It is also, with its element $\iota$,   a ``generator'' for the category of groupoids, in the similar way that the integers $\mathbb Z$ with the element $1$, form a ``generator''  for the category of groups.  Neglecting $\mathcal I$ is analogous to the long term neglect of zero in European mathematics.)

  The use of the {\it fundamental groupoid}  $\pi_1(X)$ of a space $X$, defined in terms of homotopy classes rel end points of paths $x \to y$ in $X$  was a commonplace by the 1960s. Students find it  easy to see the idea of a path as a journey, not necessarily a return journey.

 I was led to  Higgins'  paper on groupoids,  \cite{Higgins2},  for its work on free groups.  I noticed that he utilised pushouts of groupoids, and so decided to insert   in the book I was writing in 1965 an exercise on the Van Kampen Theorem for the fundamental groupoid $\pi_1(X)$. Then I  thought I had better write out a proof;  when I had done so it seemed so much better than my previous attempts that I decided to explore the relevance of  groupoids.

  It was still annoying that I  could not deduce the fundamental group of the circle!  I then realised we were in a ``Goldilocks situation'': {\it one} base point was {\it too small}; taking the {\it whole space}   was {\it too large}; but for the circle  taking {\it two}  base points was {\it  just right}! So,  we needed a definition of the fundamental groupoid $\pi_1(X,S)$ for a {\it set} $ S$ of base points chosen according to the geometry of the situation: see  the paper  \cite{BvKT}\footnote{Corollary 3.8 of that paper seems to cover the most general formula  stated in \cite{VK}, namely for a  fundamental group of the union of two spaces whose intersection is not path connected.  However Grothendieck has argued that such reductions to a group presentation may not increase understanding.}.  and  all editions of \cite{Elements}, as well as \cite{Higgins4}\footnote{For a discussion on this issue of many base points, see \url{https://mathoverflow.net/questions/40945/}.  } 
  
  The new statement of Theorem \ref{thm:vktgp} then replaces $x$ by a set $S$ meeting each path component of $W$, and we get a pushout of groupoids instead of groups. To apply this, one needs to learn how to calculate with groupoids (cf \cite{Higgins2}). 

 An inspiring  conversation with George W. Mackey (1916-2006)  in 1967 at a Swansea BMC\footnote{British Mathematical Colloquium}, where I gave an invited talk on the fundamental groupoid,  informed me of the notion of ``virtual groups'', cf \cite{Mackey,Ramsay}  and their relation to groupoids;  then led me to extensive work of C. Ehresmann (1905-1979) and his school, all showing that the idea of groupoid had much wider
  significance than I had suspected, cf \cite{Ehr05}. See also more recent work on for example  Lie groupoids, Conway groupoids, groupoids and physics.

However  the  texts on algebraic topology  which give  this   theorem  (published in \cite{BvKT})  are currently (as far as I am aware)  \cite{Elements, BHS, Z};   it is also in \cite{Higgins4}.

\section{From groupoids to higher groupoids}\label{sec;higher}
As we have shown, ``higher dimensional groups'' are just abelian groups.

However this is no longer so for ``higher dimensional groupoids'', \cite{BS1, BS2}.

It seemed to me in 1965 that some of the arguments for the VKT generalised to higher dimensions and this was prematurely claimed as a theorem in \cite{BvKT}.\footnote{It could be  more accurately  called ``There are ideas for  a proof in search of a theorem''. }

One of these arguments comes under the theme or slogan of ``algebraic inverses to subdivision''.
\begin{equation} \label{eq:aitos}
\vcenter{\xymatrix@M=0pt@=1pc{\ar @{-} [rrrrr] \ar @{-} [dddd]&&&&&\ar @{-} [dddd]\\&&&&&\\
&&&&&\\
&&&&&\\
\ar @{-} [rrrrr] &&&&&}}\qquad \leftrightarrow \qquad   \vcenter{\xymatrix@M=0pt@=1pc{\ar
@{-} [rrrrr] \ar @{-} [dddd]&\ar @{-} [dddd]&\ar @{-} [dddd]&\ar
@{-} [dddd]&\ar @{-} [dddd]&\ar @{-} [dddd]
\\\ar @{-} [rrrrr]&&&&&\\
\ar @{-} [rrrrr]&&&&&\\
\ar @{-} [rrrrr]&&&&&\\
\ar @{-} [rrrrr]&&&&&}}\end{equation}
 From left to right gives {\it subdivision}.
  From right to left should give {\it composition}.  What we need for higher dimensional, nonabelian,  local-to-global problems is:

\hspace{10em}   {\it Algebraic Inverses to Subdivision.}

This aspect is clearly more easily treated by cubical methods, rather  than the standard simplicial,  or the more recent ``globular'' ones. \footnote{Cubical subdivisions are easily expressed using a matrix notation; \cite[13.1.10] {BHS}. The ``globular'' geometry is explained in for example \cite{B-HHA}. }

One part of the proof of the VKT for the fundamental group or groupoid, namely  the uniqueness of a universal morphism,  is more easily expressed in terms of the  double groupoid $\square G$  of commutative squares in a group or groupoid $G$,  which I first saw defined in \cite{Ehresmann-65}.  The essence of its use is as follows: consider a diagram of morphisms in a groupoid:
\begin{equation}
\xybiglabels \vcenter{\xymatrix@M=0pt@=1pc{\ar @{-} [rrrrr]^b \ar @{=} [dddd]_1 &\ar @{-} [dddd]&\ar @{-} [dddd]&\ar @{-} [dddd]&\ar @{-} [dddd]& \ar @{=} [dddd] ^1 \\
\ar @{-} [rrrrr] &&&&&\\
\ar @{-} [rrrrr] &&&&&\\
\ar @{-} [rrrrr] &&&&& \\
\ar @{-} [rrrrr]_a &&&&& }}
\end{equation}
Suppose each individual square is commutative, and the two vertical outside edges are identities.
Then we easily deduce that $a = b$\footnote{The disarmingly simple higher dimensional version of this argument is \cite[13.7.5]{BHS}. }.

For the next dimension we therefore expect to need to know what is a ``commutative cube'':
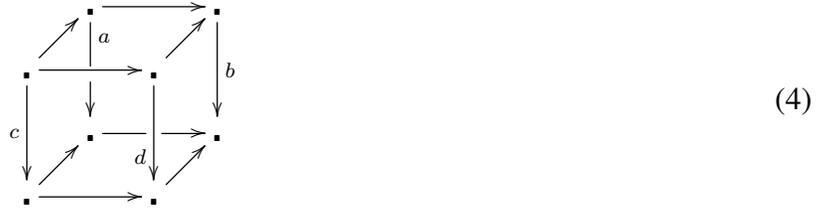
\begin{equation}\vcenter{
\xymatrix@!0{
&\bu \ar [rr]\ar ' [d]^a[dd]
& & \bu\ar [dd]^b
\\
\bu \ar [ur]\ar [rr]\ar [dd]_c
& & \bu \ar [ur]\ar [dd]_(0.68)d
\\
& \bu \ar ' [r][rr]
& & \bu
\\
\bu \ar[rr]\ar [ur]
& & \bu \ar [ur] .}}
\end{equation}
and this is expected to be in a double groupoid\footnote{More explanation is given in \cite{ B-Indag} of the way this may be given in a ``double groupoid'' of squares $G$ in which the horizontal edges $G_h$ and vertical edges $G_v$ come
from the same groupoid: i.e.  $G_h=G_v$.} We want the ``composed faces" to commute! What can this mean?

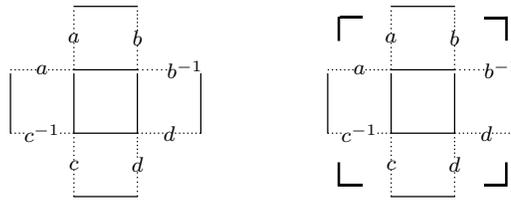
\begin{figure}[h]

$$  \xymatrix@M=0pt { & \ar@{..} [d]|a  \ar @{-} [r]  & \ar @{..} [d] |b & \\
 \ar @{-} [d] \ar @{..} [r] |a & \ar @{-}[r] \ar @{-}[d] & \ar @{-} [d] & \ar @{-} [d] \ar @{..} [l] |(0.25) {b^{-1}} \\
 \ar @{..} [r]|{c^{-1}} & \ar @{-} [r] \ar @{..} [d] |c & \ar @{..} [d] |{d}&  \ar @{..} [l] |{d} &\\
   & \ar @{-} [r] && } \qquad \xymatrix@M=0pt {\ar @{}[dr] |(0.35){\tl} & \ar@{..} [d]|a  \ar @{-} [r]  & \ar @{..} [d] |b &\ar @{}[dl]|(0.35)\tr \\
         \ar @{-} [d] \ar @{..} [r] |a & \ar @{-}[r] \ar @{-}[d] & \ar @{-} [d] & \ar @{-} [d] \ar @{..} [l] |(0.25) {b^{-1}} \\
         \ar @{..} [r]|{c^{-1}} & \ar @{-} [r] \ar @{..} [d] |c & \ar @{..} [d] |{d}&  \ar @{..} [l] |{d} &\\
         \ar @{}[ur]|(0.35)\bl & \ar @{-} [r] &&  \ar @{}[ul]|(0.35)\br }$$
\caption{``Composing" five faces of a cube}\label{fig:flatcomcub} \label{fig:commcub}
\end{figure}

 We might say that  the ``top''  face is the ``composite'' of the other
faces:  so fold them flat to give the left hand diagram of Fig. \ref{fig:flatcomcub},   where the dotted lines show adjacent edges of a ``cut". \footnote{To avoid adding an arrow to every edge of these  diagrams we adopt the convention  that edges are directed from top  to bottom and from left to right; this is why  $b^{-1}$ and $c^{-1}$ appear but not $a^{-1}, d^{-1}$.} We indicate how to glue these edges back together in the right hand diagram of this  Figure by means of extra squares which are a new  kind of ``degeneracy".

 Thus if we write the standard double groupoid identities in dimension 2 as
 $$ \tsq \quad \hh\quad \vv $$
 where a solid line indicates a constant edge, then the new types of square with commutative boundaries are written\footnote{An advantage of this ``conceptual'' or ``analogical'' notation over a more traditional or ``logical'' notation is that large diagrams involving these operations can be evaluated by eye, as in \cite[p.188]{BHS}; cf \cite{Wig} for the importance of  a conceptual approach.  However such calculations in yet higher dimensions could require appropriate computer programs!}
$$ \tl \quad \tr \quad \bl \quad \br . $$
These new kinds of ``degeneracies" were  called  {\bf connections} in  \cite{BS2}, because of a relation to path-connections in differential geometry. In a formal sense, and in all dimensions, they are constructed from  the two functions $\max, \min: \{0,1\} \to \{0,1\}$.

A basic intuition for the proof of a 2D Van  Kampen Theorem was also that  a well defined composition of  commutative cubes in any  of the three possible directions is also commutative, so this has to be proved once a full definition is set up, and then generalised to all dimensions.

It is explained in  \cite[\S 8]{B-Indag} how the use  of these ``connections'' as an extra  form of ``degeneracies'' for the traditional theory of cubical sets remedied some key  deficiencies  of the cubical as against the standard simplicial theory, deficiencies which had been known  since  1955;   the wider use of  such enhanced cubical methods then allowed better control of homotopies and higher homotopies (because of the rule $I^m \times I^n \cong I^{m+n}$). It of course kept in the cubically  allowed ``algebraic inverses to subdivision'',
and so possibilities for Higher Van Kampen Theorems; this is explained  starting in dimensions 1, 2  in \cite[Part 1]{BHS}, and continuing in Parts II, III  in all dimensions;  this book gives  what amounts to a rewrite of much traditional singular simplicial  and cellular algebraic topology, cf \cite{B-Indag}.

\section{The influence of work 1941-1949 of J.H.C. Whitehead}
The relation of the  fundamental group to aspects of geometric group theory was an important feature of the work of Poincar{\' e}. The relation of various versions of homotopy groups to group theory was an important feature of the work of  J. Henry C. Whitehead (1904-1960), my supervisor in the period 1957-1959.  His paper  \cite{CHI}  is known as basic in homotopy theory, and the use of the word ``combinatorial" in its title indicates its links with combinatorial group theory. His paper \cite{CHII} is less well known, but is the basis of \cite{ML-W}, which describes homotopy 3-types (now called 2-types) in terms of the algebra of ``crossed modules'', \eqref{sec:xmod}.

I overheard Whitehead tell J.W. Milnor (1931- ) in 1958 that the early homotopy theorists were fascinated by the operations of the fundamental group on the higher homotopy groups, and also by the problem of computing the latter, preferably with this action.

Whitehead was able by very hard work and study to look at an area and seek out major problems.   One of  his aims in the late 1930s was to discover whether the Tietze transformations of combinatorial group theory could be ``extended''  to higher dimensions. His main method of such extension was envisaged as ``expansions'' and ``collapses'' in a simplicial complex;  he also wrestled with the problem of simplifying such complexes into some kind of ``membrane complex''; in his work after the war this became codified in \cite{CHI} as the notion of ``CW-complex'', and his work on generalising Tietze transformations became a key part of algebraic $K$--theory.

He was very concerned with the work of Reidemeister
   on the relations  between simplicial complexes and presentations of groups, and methods of finding appropriate geometric  models of group constructions, particularly generators and relations, and possible higher  analogues. It was only after the war that the topological notion of ``adding a 2-cell'', as compared to ``adding a relation'',  was gradually codified  through  the notion of ``adjunction space'',   \cite{W-B}. This   gives  a useful  method   of constructing a space $Y$ as an identification space $B \cup_f X$  of $X \sqcup B$ given  a space  $X$  with an inclusion $ i: A \to X$, and a map $f: A \to B$, yielding a $g: X \to Y$. This definition, which allows for constructing continuous functions from  $Y$, was background to   the notion of CW-complex in \cite{CHI}. A basic account of adjunction spaces and their use in homotopy theory  is in \cite{Elements}.

It was gradually realised, cf \cite{Whi46},   that for a pair $(X,A)$ of pointed spaces, i.e. a space $X$,  subspace $A$ of $X$,  and point $x \in A$, there was an exact sequence of groups which ended with
\begin{equation}\label{eq:2rel}
  \pi_2(X, x) \to \pi_2(X,A,x) \xrightarrow{\delta} \pi_1(A,x) \xrightarrow{i} \pi_1(X,x) \to 0
\end{equation}
where $\pi_2(X,A,x) $ is the {\it second  relative homotopy group} defined as homotopy classes rel vertices of maps $(I^2,J,V)\to (X,A,x)$ where $V$ are the vertices  of the unit square $I^2$, and  $ J$  consists of all edges except one, say $\partial^-_1 I^2$; with this choice, the composition of such classes is taken in direction $2$, as in:

\begin{equation}\label{eq:relhomgps}\vcenter{\xybiglabels
\xymatrix@M=0pc{\ar @{=} [rr]^x \ar @{=}[d]_x \ar @{} [dr]|X&\ar@{=}[d]  \ar @{..}[r] \ar @{} [dr]|X &\ar @{=}[d] \ar @{--}[rr]^x   \ar @{} [drr]|X &   & \ar @{} [dr]|X\ar  @{=}[d] \ar  @{=} [rr]^ x  & \ar @{} [dr]|X   \ar  @{=}[d]  &     \ar  @{=}[d]^x \\
\ar @{}[r]|A&\ar @{}[r]|A &\ar @{}[rr] |A& &\ar @{}[r]|A&\ar @{}[r]|A&}}\quad \sdirects{2}{1}
  \end{equation}

The groups on either side of $\xrightarrow{\delta}$  in the sequence  \eqref{eq:2rel} are in general nonabelian. Whitehead saw that there was an operation $(m,p)\mapsto m^p$ of $P=\pi_1(A,x) $ on $M=\pi_2(X,A,x) $ such that $\delta(m^p)= p^{-1}(\delta m)p$ for all $ p \in P, m \in M$. In a footnote of  \cite[ p.422 ]{Whi41a} he also stated   a rule equivalent to $m^{-1}nm = n^{\delta m }, m,n  \in M$. The standard proof of this rule uses a ``2-dimensional argument''  which can be shown in the following diagram in which $a=\delta m, b = \delta n$:

\begin{equation}\label{eq:relhomgps}
\xybiglabels
  \vcenter{\xymatrix@M=0pt@=3pc{\ar @{=}[r] \ar @{=} [dd]\ar @{}[dr]|{-_2m} &\ar @{=} [dd] \ar @{=} [r]{} \ar @{}[dr]|\sq &\ar @{=} [dd] \ar @{}[dr]|{m}\ar @{=} [r] & \ar @{=}[dd]\\
  \ar @{-} [r] \ar @{}[dr]|{\vv} &\ar @{=} [r]\ar @{} [dr]|{n} &\ar @{-} [r] \ar @{}[dr]|{\vv}& \\   \ar @{-} [r]_{a^{-1}}  &\ar @{-}[r]_b&\ar [r]_a & }} \qquad \sdirects{2}{1}
\end{equation}
Here the double lines indicate constant paths, while $\vv, \sq$ denote a vertical and double identity respectively. He later  introduced the term {\it crossed module} for this structure, which has  an important place in our story.

Whitehead realised that to calculate $\pi_2(X,x)$ we were really in the business of calculating the group morphism $\delta$ of \eqref{eq:2rel}, and that such a calculation should involve the crossed module structure. So one of his strengths, as I see it, was that he was always on the lookout  for the controlling {\it underlying structure}. In particular, \cite[\S 16]{CHII} sets up the notion of {\it free crossed module} and uses geometric methods from  \cite{Whi41a} to show how such can be obtained by attaching cells to a space. This theorem on free crossed modules, which  is sometimes quoted but seldom proved in even advanced texts on algebraic topology, was a direct stimulus to the work on higher homotopy groupoids.

 As explained in \cite[\S 8]{B-Indag}, Higgins and I agreed  in 1974  that this result on free crossed modules was a rare, if not the only, example of a universal nonabelian property in 2-dimensional homotopy theory. So if our conjectured but not yet  formulated  theory was to be any good, it should have Whitehead's theorem as a Corollary. But that theorem was about second {\it relative} homotopy groups. Therefore we  also should look at a relative situation, say  $S \subseteq A \subseteq X$ where $S$ is a set of base points and $X$ is a space.

There is then a simple way of getting what looks like a putative homotopy double groupoid from this situation:    instead of the necessarily ``1-dimensional composition'' indicated in diagram \eqref{eq:relhomgps},   we should take the unit square $I^2= I \times I$,  where $I$ is the unit interval $[0.1]$, with  $E$ as its set of edges, $V$ as its set of vertices;  we should  consider the set $R_2(X,A,S)$ of maps $(I^2, E, V) \to (X,A,S)$, and then take homotopy  classes relative to the vertices $V$ of such maps to form say $\rho_2(X,A,S)$.  It was this last set that was now fairly easily shown to have the structure of double groupoid over $\pi_1(A,S)$ {\it with connections}, and so to be  a 2-dimensional version of the fundamental groupoid on a set of base points!

I need to explain the term {\it connections}, as introduced in dimension 2
 in \cite{BS2}. It arose from the desire  to construct examples of double groupoids  other than the  previously defined $\square G$  of commutative squares in a groupoid $G$.

Whitehead  proved in \cite{Whi46} that the boundary $\delta:  \pi_2(X,A, x) \to \pi_1(A,x)$ and an action of the group $\pi_1(A,x)$ on the group $\pi_2(X,A,x)$ has the structure of  a crossed module. He also proved in \cite[\S 16]{CHII} what we call Whitehead's {\it free crossed module theorem } that in the case $X$ is formed from $A$ by attaching $2$-cells, then this crossed module is {\it free} on the characteristic maps of the attaching 2-cells;  this  topological model of ``adding relations to a group'' is sometimes stated but rarely proved in texts on algebraic topology. Later, an {\it exposition}  of the proof as written out  in \cite[\S 16]{CHII} was published as \cite{B-80}; it uses  methods of  knot theory  (Wirtinger presentation) developed in the 1930s and of ``transversality'' (developed  further in  the 1960s). but taken from the    papers \cite{Whi41b,Whi46}. The result was earlier put in the far wider  context of a 2-dimensional Van Kampen type Theorem in  \cite{BH78}.

The notion of crossed module occurs in other algebraic contexts: cf \cite{Lue}, which refers also to 1962 work of the algebraic number theorist A. Fr{\" o}hlich on nonabelian homological algebra, and more recently in  for example \cite{Jan,MVL}.

 The  definition
  of this crossed module  in \eqref{eq:relhomgps} involves choosing which vertex should be  the base point of the square and which edges of the square should map to the base point $a$, so that the remaining edge maps into $A$. However it  is a good principle  to reduce, preferably completely, the number of choices used in basic definitions (though such choices are likely in developing consequences of the definitions).  The paper \cite{BH78}, submitted in 1975, defined
  for a triple of spaces $\mathbf X= (X,A,S)$ of spaces such that $S \subseteq A \subseteq X$  a structure $\rho(\mathbf X )$: this  consisted  in dimension $0$ of $S$ (as a set);  in dimension $1$ of $ \pi_1(A,S)$;  and in dimension $2$ of homotopy classes relative  to the vertices of maps $(I^2,   E, V) \to (X, A, S)$, where $E,V$ are the spaces of edges and vertices of the standard square.\footnote{In that paper  it was assumed that each loop in $S$ is contractible in $A$ but this later proved too restrictive on $A$,  and so  homotopies fixed on the vertices of $I^2$ and in general $I^n$, were used in \cite{BHS}.}

   \begin{equation}\label{eq:htydb}\xybiglabels
     \vcenter{\xymatrix{S \ar @{-} [r] | A  \ar @{}[dr]|X \ar @{-}[d] |A & S\ar@{-} [d] |A \\S \ar @{-} [r] | A & S}    } \qquad \sdirects{2}{1}
      \end{equation}

  This definition makes {\it no choice of preferred direction}. It is fairly  easy and direct to prove that $\rho(\mathbf X)$ may be given the structure of double groupoid with connection\footnote{Or alternatively, of double grouoid with thin structure,  \cite[p.163]{BHS}. } containing  a copy of  the double groupoid $\sq\pi_1(A,S)$. That is, the proofs of the required properties of $\rho(\mathbf X) $ to make it a $2$-dimensional  version of the fundamental group  as sought in the 1930s are fairly easy but not entirely trivial. The longer task, 1965-1974,  was formulating the ``correct'' concepts (in the face of prejudice from some referees and editors). The proof of the  corresponding Van Kampen Theorem allows a  nonabelian result in dimension 2 which vastly  generalises the work  of \cite[\S 16]{CHII}: for example, it gives a result when $X$ is formed from $A$ by attaching a cone on $B$, Whitehead's case being when $B$ is a wedge of circles.  See also \cite[Chapter 5]{BHS} for many other explicit homotopical excision examples, some using the GAP system.

  Here is an example where we use more than one base point. Let $X$ be the space $S^2 \vee [0,1] $ where $0$ is identified with $N$, the North pole of the 2-sphere. Let $S= \{ 0,1 \}$ as  a subset of $X$. We know that $\pi_2(X,0) \cong \mathbb Z$ and it is easy to deduce that $\pi_2(X, S)$ is the free $\mathcal I =  \pi_1([0,1], S)$-module on one generator. Now we can use the relevant 2-D Van Kampen Theorem to show that $\pi_2(S^2 \vee S^1, N)$ is the free $\mathbb Z$-module on one generator.

Note also that it is easy to think of generalisations to higher dimensions of diagram \eqref{eq:htydb}, namely to the filtered spaces of \cite{BHS}, following the lead of \cite{Bl48}\footnote{That paper, as does \cite{CHII},  uses the term ``group system '' for what is later called a ``reduced crossed complex'', i.e. one with a single base point.  }.

\section{General considerations }
The paper \cite[\S 2]{B-Indag}  argues that one difficulty of obtaining such strict higher structures, and so  theorems on colimits  rather than {\it homotopy}  colimits,  is the difficulty of working with ``bare''  topological spaces, that is topological spaces with no other structure (the term ``bare'' comes from \cite[\S 5]{GrEsq}).

The argument is the practical one that in order to calculate an invariant of a space one needs some information on that space:  that information will have a particular  algebraic or geometric structure which needs to be used. Because of the variety of convex sets in dimensions higher than $1$,  there is a variety also of potentially relevant higher algebraic structures. It turns out that some of these structures are non trivially equivalent, and can be described as ``broad'' and ``narrow''.

The broad ones are elegant and symmetric, and useful for conjecturing and proving theorems; the narrow ones are useful for calculating and relating to classical methods;  the  non trivial equivalence allows one to get the best use of both. An example in \cite{BHS} is the treatment of cellular methods, using filtered spaces and the related algebraic structures  of crossed complexes (``narrow''), and cubical   $\omega$-groupoids (``broad'').

This use of structured spaces is one explanation of why the account in \cite{BHS} can, in the tradition of homotopy theory,  use and calculate with,  strict rather than lax, i.e. up to homotopy, algebraic structures. In comparison, the paper \cite{BKP} does give a strict result for all Hausdorff spaces, but so far has given no useful consequences.

There is a raft of other papers following Grothendieck in exploiting the idea of  using ``lax structures'',  involving homotopies,   homotopies of homotopies, $\ldots$,  in a simplicial context; this involves seeing the simplicial singular complex $S^\Delta(X)$ of a space $X$ as a form of ``$\infty$-groupoid'', and which have had some famous applications.

 To get nearer to a fully nonabelian theory we so far have only the use of $n$-cubes of pointed spaces as in \cite{BL1,BL2,JFA}. It is this restriction to {\it pointed spaces}  that is a kind of anomaly, and has been strongly criticised by Grothendieck as not suitable for modelling in algebraic geometry. However  the paper \cite{ESt} gives an application to a well known problem in homotopy theory, namely  {\it  determining   the first non-vanishing homotopy group of an $n$-ad}; also   the {\it nonabelian tensor product of groups}  from \cite{BL1} has become  a flourishing topic in group theory (and analogously for Lie algebras);   a bibliography\footnote{See \url{http://www.groupoids.org.uk/nonabtens.html}} 1952-2009 has 175 items.

  The title of the paper \cite{hdgt} was also  intended to stimulate the  intuition that ``higher dimensional geometry requires higher dimensional algebra'', and so to encourage  non rigid argument on the forms that the latter  could and should take. Perhaps the early  seminar of Einstein \cite{Ein22} could be helpful in this.

   It is now a commonplace that the further development of related  higher structures are important for mathematics and particularly for applications in physics\footnote{This assertion is supported  by a web search on ``Institute of higher structures in maths''. }. Note that the mathematical notion of group is deemed fundamental to the idea of symmetry, whose implications range far and wide. The bijections of a set $S$ form a group $Aut(S)$. The automorphisms $Aut(G)$ of a group $G$ form part of a crossed module $\chi: G \to Aut(G)$. The automorphisms of a crossed module form part of a ``crossed square'' \cite{BG}. These structures of set, group, crossed module, crossed square, are related to homotopy $n$-types for $n=0,1,2,3$.

The use  in texts on  algebraic topology  of {\it sets of  base points}  for fundamental groupoids seems currently restricted to \cite{Elements,BHS,Z}.

The argument over {\v C}ech's seminar to the 1932 ICM seems now able to be resolved through this development of groupoid and higher groupoid work, and he  surely deserves  credit for the first presentation on higher homotopy groups, as reported in  \cite{Alex, Ce,jam}.

Another way of putting our initial quotation from Einstein is that one should be wary of ``received wisdom''.

\section{Appendix}
 In this section we give a glimpse of some of the calculations needed to show the axioms of a crossed module do work to give rise to a double groupoid.
\subsection{Crossed modules and  compositions of labled squares }\label{sec:xmod}
An  easy  example of a double groupoid is to start with a group $P$ and consider the set $\sq P$ of commuting squares in $P$, i.e. quadruples $(a  ^c _ b d)$ such that $ab=cd$.   Any well defined composition of commuting squares is commutative; for example  $ab=cd $ and $dg=ef $ implies $ abg= cef$ . So it is easy to see that $ \sq P$ forms the structure of double groupoid.

In homotopy theory we do not expect all squares of morphisms to commute. So it is sensible to consider a  subgroup say $M$ of $ P$  and squares which commute up to an element of $M$. There are many choices to make here; suppose we make the convention that we consider squares
\begin{equation}\xybiglabels
\vcenter{  \xymatrix@M=0pc{
  \ar [r] ^g \ar[d] _h \ar @{} [dr] |m & \ar [d] ^a \\
  \ar [r] _k &\cdot }}  \quad \xdirects {2}{1}
\end{equation}
in which $a,g,h,k \in P, m \in M$ and $k^{-1}h^{-1} ga=m $. So we are starting with the bottom right hand corner as `base point',  and going clockwise around the square. You quickly find that for a  composition of such squares to work you need the subgroup $M$ to be normal in $P$.

In homotopy theory you expect many ways of making a boundary commute. So it seems sensible to replace a subgroup $M$ of $P$ by a morphism $\mu: M \to P$. It  also seems sensible to replace the group $P$ by a groupoid. What then should be the conditions on $\mu$? Convenient ones turned out to be a groupoid version of a notion  envisaged in a footnote of the paper \cite[p.422]{Whi41a};  in \cite{Whi46}  this structure was called a {\it crossed module} and it was further  developed in \cite{CHII}.
\begin{Def} A morphism  $\mu: M \to P$   groupoids is called a  {\it crossed module} if $\mu$ is the identity on objects;  $M$ is discrete, i.e. $M(x,y) $ is empty for $x \ne y$;  $P$ operates on  the right  of the group $M$, $(m,p) \mapsto m^p$. satisfying  the additional rules to those of an operation:

CM1)  $\mu (m^p) = p^{-1} \mu (m) p$;

CM2) $n^{-1} m n= m ^{\mu n}$\\
for all $m,n \in M, p \in P$.
\end{Def}

We start with a crossed module $\mu: M \to P$ where $P$ is a groupoid with object set $S$ and $M$ is a discrete groupoid consisting of groups $M(s), s \in S$ on which $P$ operates.

We then form the elements of our double groupoid to be in dimension $0$, the elements of $S$; in dimension $1$ the elements of $P$; and in dimension $2$ the quintuples consisting of  one element of $M$ and four elements of $P $  whose geometry forms a square as in the left hand diagram below and $\mu(n ) = k^{-1}h^{-1} ga$.

Further we define  such a ``filled square'' to be {\it thin} if $n= 1$. Thus the thin elements form a special kind of commutative square.

 We try a `horizontal' composition
\begin{equation}\xybiglabels
\vcenter{  \xymatrix@M=0pc{
  \ar  [r] ^g \ar[d] _h \ar @{}[dr] |n & \ar  [d] ^a \\
  \ar [r] _k &\cdot }} \circ_2 \vcenter{\xymatrix@M=0pc{
  \ar  [r] ^c  \ar [d] _a \ar @{}[dr] |m & \ar  [d] ^d \\
  \ar  [r] _b &\cdot }} = \vcenter{  \xymatrix@M=0pc{
  \ar  [r] ^{gc} \ar[d] _h \ar @{}[dr] |\alpha & \ar  [d] ^d \\
  \ar  [r] _{kb}  &\cdot }} \qquad \sdirects{2}{1}
 \end{equation} assuming $gc, kb$ are defined,
  and a `vertical' composition
  \begin{equation}\xybiglabels
\vcenter{  \xymatrix@M=0pc{
  \ar [r] ^g \ar [d] _f \ar @{}[dr] |u& \ar  [d] ^e \\
  \ar [r] _c &\cdot }} \circ_1 \vcenter{\xymatrix@M=0pc{
  \ar  [r] ^c \ar [d] _a \ar @{}[dr] |m & \ar [d] ^d \\
  \ar  [r] _h &\cdot }} = \vcenter{  \xymatrix@M=0pc{
  \ar [r] ^{g} \ar [d] _{fa} \ar @{}[dr] |\beta  & \ar  [d] ^{ed} \\
  \ar  [r] _{h}  &\cdot }} \qquad \sdirects{2}{1}
 \end{equation}
 assuming $fa, cd$ are definied.
 The problem is to give values for $\alpha, \beta $   and to  prove in each case that the square fits the definition. In fact we find that $\alpha =(n^b) m$ and $\beta = mu^d$ will do the trick; these calculations  strongly use the two rules for crossed modules.

 I won't give the argument here, as it is not hard, and quite fun, and  is on pages 176-178 of \cite{BHS}, including a full proof, which again needs both axioms CM1) and CM2),  of the interchange law for $\circ_1, \circ_2$.

  The thin elements are related to  a functor from the category $ \mathbf{CM}$ of  crossed modules  to the category $\mathbf{DGT}$ of double groupoids with thin structure, which is an equivalence of categories, \cite{BS2,Higgins05}

  These two categories play different roles, and  in the language of \cite{B-Indag} they are  called {\it narrow } and {\it broad} algebraic structures respectively; the `narrow' category is used for calculation and relation to traditional terms  such as, in our case,  homotopy groups. The `broad' category is used for expressive work, such as  formulating ``higher diimensional composition,  conjectures and proofs. The equivalence between the two categories, which may  entail a number of somewhat arbitrary choices,  enables us  to use whichever is convenient for the job at hand, often without worrying about the details of the proof of the equivalence. This is especially important in the case of dimensions $> 2$.  Such a  use of equivalent categories to provide different types of tools for research should perhaps be considered as part of the ``methodology'' of mathematics'', \cite{BP}.

  The discussion on \cite[p.163]{BHS} relates crossed modules to ``double groupoids with thin structure''; the category of those has  the advantage of being ``Yoneda invariant'' and so can be repeated internally in any category with finite limits. Wider applications of the concept of crossed module are also shown in for example \cite{MVL}, and of double groupoid with thin structure in  \cite{FMP}. Perhaps this distinction between ``narrow'' and ``broad'', which becomes even more stark  in higher dimensions,   is relevant to  Einstein's old discussion  in \cite{Ein22}.

  I think it is fair to say that the higher SVKT's would not have been formulated, let alone proved,  in a ``narrow'' category, cf \cite{B-Indag},  while the formulation in the novel ``broad'' category was initially felt by one editor to be ``an embarrassment''\footnote{We refer here to the comment in  \cite[p. 48]{Rota} in reply  to the question: "What can you prove with exterior algebra that you cannot prove without it?".  Rota  retorts:   `` Exterior algebra is not meant to prove old facts, it is meant to disclose a new world. Disclosing new worlds is as worthwhile a mathematical enterprise as proving old conjectures. ''}.  The fact that an SVKT is essentially a colimit theorem implies of course that it should  give precise algebraic calculations, not obtainable by means of say exact sequences or spectral sequences. The connectivity conditions for the use of the SVKT also limit its applicability;  the 2-dimensional theorem does enable some computations of homotopy 2-types; but as has been known since the 1940s in terms of group theory considerations, computation of a morphism does not necessarily imply computation of its kernel.

   Thus the search started by the early topologists has a resolution:  precise higher dimensional versions of the Poincar\' e fundamental group, using higher groupoids, do exist in all dimensions, as described in \cite{BHS}.  Whitehead's work on free crossed modules is there seen, \cite[p.235]{BHS}, as a special case of a result on  ``inducing'' a new  crossed module $ f_*M \to Q$ from a crossed module  $\mu : M \to P$ by a morphism $f: P \to Q$ of groups (or groupoids), and this result itself is a special case of a Van Kampen type theorem, which includes classical theorems such as the Relative Hurewicz Theorem.  The important point is that the theory allows caculations and applicati9ons not previously possible,

   Those particular structures however model  only a limited range of homotopy types. There is another theory for pointed spaces  due to Loday, which is proved in \cite{BL1} also to have a higher SVKT. It has yielded a range of new applications, \cite{ESt,JFA},  but its current limitation   to pointed spaces makes it less suitable for other  areas, such as algebraic geometry.

   Note that Whitehead's paper  \cite{SHT} is a sequel to \cite{CHI,CHII} and other earlier papers. but the treatment is different.

   There is also a large literature on models using ``lax  higher homotopy groupoids'', stimulated by ideas of Grothendieck, see \cite[p.xiv] {BHS}.

\section*{Acknowledgements}
The projects described here have been supported by grants from SRC, EPSRC, British Council, University College of North Wales, Strasbourg University, Royal Society, Intas.

The perceptive comments of two referees have been helpful in developing this paper.

I am grateful to Terry Wall who showed me in 1975 how our cubical  method  of   constructing  a well defined strict homotopy groupoid in dimension 2 might be extended to at least one higher dimension by modelling  ideas from \cite{Whi41b},  leading eventually to the collapsings and expansions of chains of cubical partial boxes of \cite[p. 380]{BHS}. The acknowledgements given in \cite{BHS} also apply of course to this paper; I again thank  Chris D. Wensley for detailed comments.

 \newcommand{\enquote}[1]{`#1'}

\end{document}